\documentclass[ECP]{ejpecp} 

\usepackage{enumerate}  
\usepackage{color}
\theoremstyle{definition}
\allowdisplaybreaks

\SHORTTITLE{2d delta Bose gas in a weighted space}

\TITLE{Two-dimensional delta Bose gas in a weighted space\support{%
	SS is partially supported by the NSF through DMS-2243112.
	LCT is partially supported by the NSF through DMS-2243112 and the Alfred P.\ Sloan Foundation through the Sloan Research Fellowship FG-2022-19308.%
}
}

\AUTHORS{%
	Sudheesh~Surendranath\footnote{University of Utah
		\EMAIL{sudheesh@math.utah.edu}}
	\and 
	Li-Cheng~Tsai\footnote{University of Utah \BEMAIL{licheng.tsai@utah.edu} }}


\AMSSUBJ{%
46N30  	
} 

\SUBMITTED{November 5, 2024} 
\ACCEPTED{April 8, 2025} 

\VOLUME{30}
\YEAR{2025}
\PAPERNUM{35}
\DOI{10.1214/25-ECP685}

\ABSTRACT{%
We extend the construction of the semigroup of the two-dimensional delta-Bose gas in \cite{gu2021moments} (based on \cite{rajeev99,dimock04}) to a weighted $L^2$ space that allows exponentially growing functions.
We further show that the semigroup of the mollified delta-Bose gas converges strongly to that of the delta-Bose gas.
}

\newcommand{\e}{\varepsilon}

\newcommand{\vectau}{\vec{\tau}}

\newcommand{\heatsg}{\mathcal{P}}	
\newcommand{\sg}{\mathcal{Q}}		
\newcommand{\ham}{\mathcal{H}}		
\newcommand{\sgsum}{\mathcal{R}}	
\newcommand{\Jop}{\mathcal{J}}	
\newcommand{\jfn}{\mathsf{j}}
\newcommand{\weiop}{\mathcal{E}}	
\newcommand{\hk}{\mathsf{p}}		
\newcommand{\Top}{\mathcal{T}}		
\newcommand{\Bsp}{\mathcal{B}}		

\newcommand{\pair}{\mathrm{Pair}}
\newcommand{\intv}[1]{[\![#1]\!]}
\newcommand{\dgm}{\mathrm{Dgm}}
\newcommand{\vecalpha}{\vec{\alpha}}

\newcommand{\veceta}{\vec{\eta}}
\newcommand{\textc}{\mathrm{c}}
\newcommand{\textr}{\mathrm{r}}

\newcommand{\yc}{y_\mathrm{c}}

\newcommand{\yr}{y_\mathrm{r}}

\newcommand{\Fou}{\widehat}				
\newcommand{\pFou}{\overbracket[.5pt]}	

\newcommand{\norm}[1]{\Vert #1\Vert}
\newcommand{\Norm}[1]{\big\Vert #1\big\Vert}
\newcommand{\NOrm}[1]{\Big\Vert #1\Big\Vert}
\newcommand{\norma}[1]{\Vert #1\Vert_{2,a}}
\newcommand{\Norma}[1]{\big\Vert #1\big\Vert_{2,a}}

\newcommand{\normop}[1]{\Vert #1\Vert_{2\to 2}}
\newcommand{\Normop}[1]{\big\Vert #1\big\Vert_{2 \to 2}}
\newcommand{\NOrmop}[1]{\Big\Vert #1\Big\Vert_{2 \to 2}}
\newcommand{\normopa}[1]{\Vert #1\Vert_{2,a\to 2,a}}
\newcommand{\Normopa}[1]{\big\Vert #1\big\Vert_{2,a\to 2,a}}
\newcommand{\NOrmopa}[1]{\Big\Vert #1\Big\Vert_{2,a\to 2,a}}
\newcommand{\ip}[1]{\langle #1\rangle}
\newcommand{\Ip}[1]{\big\langle #1\big\rangle}

\newcommand{\normone}[1]{\vert #1 \vert_{1}}
\newcommand{\normtwo}[1]{\vert #1 \vert_{2}}

\newcommand{\C}{\mathbb{C}}

\newcommand{\R}{\mathbb{R}}
\newcommand{\Z}{\mathbb{Z}}
\newcommand{\Csp}{C}
\newcommand{\Lsp}{L}

\newcommand{\ind}{\mathbf{1}}			
		
\renewcommand{\d}{\mathrm{d}}		
\renewcommand{\bar}{\overline}

\newcommand{\til}{\widetilde}
\newcommand{\img}{\mathbf{i}}
\renewcommand{\Re}{\operatorname{Re}}
\renewcommand{\Im}{\operatorname{Im}}

\begin{document}
	

	
\section{Introduction}
\label{s.intro}
In this paper, we study the two-dimensional, $n$-particle delta-Bose gas.
Fix a nonnegative $\Phi\in\Csp^\infty(\R^2)$ with a compact support such that $\int_{\R^2}\d x\,\Phi(x)=1$, and, for $\e\to 0$, let $\Phi^\e(x):=\e^{-2}\Phi(\e^{-1}x)$ be the corresponding approximation of the delta function.
Let $\intv{n}:=\{1,\ldots,n\}$, write $\alpha=ij=\{i<j\}$ for an unordered pair of integers in $\intv{n}$, view $ij$ as a \emph{set}, and let $\pair\intv{n}:=\{ij|ij\subset\intv{n}\}$.
Consider $\R^{2\intv{n}}:=(\R^2)^{\intv{n}}=\{x=(x_1,\ldots,x_n)| x_i\in\R^2\}$, and, for $\alpha=ij$, let $\Phi^\e_\alpha(x):=\Phi^{\e}(x_i-x_j)$.
The Hamiltonian $\ham^{\e}$ of the mollified delta-Bose gas acts on functions on $\R^{2\intv{n}}$ and is given by
\begin{align}
	\label{e.ham}
	-\ham^{\e} :=  \frac{1}{2}\sum_{i=1}^n \Delta_{i} + \beta_\e \sum_{\alpha\in\pair\intv{n}}\Phi^{\e}_{\alpha} \ ,
\end{align}
where $\Delta_i$ is the Laplacian acting on $x_i\in\R^2$, $\Phi^{\e}_{\alpha}$ acts multiplicatively,  
\begin{align}
	\label{e.beta}
	\beta_\e
	&:=
	\frac{2\pi}{|\log\e|}
	+
	\frac{\pi}{|\log\e|^2} \Big(
		\theta - 2\log 2 +2\gamma +2 \int_{\R^4} \d x \, \d x' \Phi(x) \log |x-x'| \Phi(x')
	\Big)\ ,
\end{align}
$\gamma=0.577\ldots$ denotes the Euler–Mascheroni constant, and $\theta$ is a fixed parameter that can be taken to be any real value.
In \eqref{e.beta}, the constants after $\theta$ are just conventional (and can be absorbed into $\theta$).

The delta-Bose gas is an instance of quantum many-body systems and is highly relevant in the study of the Stochastic Heat Equation (SHE).
In the quantum context, \eqref{e.ham} can be used to approximate the Hamiltonian of $n$ Bosons with pairwise, attractive (true) delta potentials.
Quantum systems with delta potentials exhibit intriguing physical properties, pose significant mathematical challenges, and enjoy a long history of study.
We point to \cite{albeverio88} for a review and mention the works \cite{dellantonio1994hamiltonians, rajeev99, dimock04, griesemer2020short, griesemer2022short} on the two-dimensional delta-Bose gas.
As it turns out, in two dimensions, in order to see a non-degenerate limit as $\e\to 0$, one needs to attenuate the coupling constant $\beta_\e$ logarithmically around the value $2\pi$ as in \eqref{e.beta}.
This was observed in the $n=2$ case in \cite[Chapter~I.5]{albeverio88} and \cite{bertini98}.
In the SHE context, the semigroup $\sg^{\e}(t):=e^{-t\ham^{\e}}$ describes the moments of the noise-mollified SHE.
Based on \cite{albeverio88,albeverio1995fundamental}, the work \cite{bertini98} established the convergence (as $\e\to 0$) of the second moment of the noise-mollified SHE to an explicit limit.
Convergence of higher moments are significantly more challenging, and was later obtained in \cite{caravenna18} for the third moment based on \cite{caravenna18b}, and in \cite{gu2021moments} for of all moments based on \cite{rajeev99,dimock04}; see \cite{chen2022two, chen2024delta, chen2024stochastic} for probabilistic studies of the moments. 
A central question in this context concerns the $\e\to 0$ limit of the noise-mollified SHE, and the limit is called the Stochastic Heat Flow (SHF).
The construction of the SHF had been a challenging open problem and was accomplished \cite{caravenna2023critical}. 

In this paper, we generalize the results of \cite{gu2021moments} from the $\Lsp^2$ space to a \emph{weighted} $\Lsp^2$ space.
Based on \cite{rajeev99,dimock04}, the work \cite{gu2021moments} gave an explicit formula of a strongly continuous semigroup $\sg(t)$ on $\Lsp^2(\R^{2\intv{n}})$ and showed that $\sg^{\e}(t):=e^{-t\ham^{\e}}$ converges strongly to $\sg(t)$ on $\Lsp^2(\R^{2\intv{n}})$.
These results translate into the convergence of the moments of the noise-mollified SHE to an explicit limit, for $\Lsp^2$ initial conditions.
The class of $\Lsp^2$ functions is natural, because the generator of $\sg(t)$ is self-adjoint, but the class leaves out some commonly considered initial conditions, most notably the constant/flat ones.
Further, when using delta-Bose gas to study the noise-mollified SHE and the SHF, one often needs to consider more general spaces than $\Lsp^2$.
One instance is the work \cite{caravenna2023critical}, which introduces weighted $\Lsp^p$ estimates for the semigroup in a discrete setting. Another instance is the recent work \cite{tsai2024}, which requires the use of the weighted $\Lsp^2$ space and the bounds \eqref{e.comp.e}--\eqref{e.bds.e} from this paper.

To state our result, let $|\cdot|_{p}$ denote the $\ell^p$ norm on $\R^{d}$, for $a\in\R$, let
\begin{align}
	\label{e.norma}
	\norma{f}
	:=
	\Big( \int_{\R^{2\intv{n}}} \d x \, \big| f(x) e^{a\normone{x}} \big|^2 \Big)^{1/2} \ ,
	\qquad
	\Lsp^{2}_{a}(\R^{2\intv{n}})
	:=
	\big\{ f \, \big|\, \norma{f}<\infty \big\}\ ,
\end{align}
let $\ip{f_1,f_2}:=\int_{\R^{2\intv{n}}}\d x\, e^{2a\normone{x}}\bar{f_1(x)} f_2(x)$ denote the inner product, and for a linear operator $\Top:\Lsp^{2}_a(\R^{2\intv{n}})\to\Lsp^{2}_a(\R^{2\intv{n}})$, let
\begin{align}
	\label{e.normopa}
	\normopa{\Top}
	:=
	\sup\big\{ \ip{f',\Top f} \, \big| \, \norma{f'} \leq 1 \, , \norma{f}\leq 1 \big\}
\end{align}
denote the weighted operator norm.
Deferring the definition of $\sg(t)$ to Section~\ref{s.limiting}, we now state our main result.
\begin{theorem}\label{t.main}
Given any $\theta\in\R$, $n\in\Z_{>0}$, and $a\in\R$, there exists $c=c(\theta,n,a,\Phi)$ such that for all $t \geq 0$, $\e\in(0,1/c]$, $f\in\Lsp^{2}_{a}(\R^{2\intv{n}})$, and $T<\infty$,
\begin{enumerate}
\item \label{t.main.1}
$\normopa{\sg(t)} \leq c\,e^{ct}$,
\item \label{t.main.2}
$\normopa{\sg^{\e}(t)} \leq c\,e^{ct}$, and $\norma{\sg^{\e}(t)f-\sg(t)f} \to 0$ as $\e\to 0$ uniformly over $t\in[0,T]$.
\end{enumerate}
\end{theorem} 
\noindent{}%
We note in the passing that the work \cite{chen2024delta} studied a similar problem for bounded $f$s and proved the pointwise (in $x$) convergence of $(\sg^{\e}(t)f)(x)$.

A crucial property used in \cite{gu2021moments} is that the generators of $\sg(t)$ and $\sg^{\e}(t)$ are self-adjoint in $\Lsp^2$, but this property no longer holds in $\Lsp^2_{a}$.
In the self-adjoint setting of \cite{gu2021moments}, to prove statements like those in Theorem~\ref{t.main}, one can just work with the resolvents. 
Indeed, in the self-adjoint setting, thanks to the spectral theorem, bounds on the resolvents imply bounds on the spectra and hence the semigroups.
In our non-self-adjoint setting, we need to bound the semigroups  $\sg(t)$ and $\sg^\e(t)$  directly, without resorting to resolvents.
These semigroups are much less well-behaved than their resolvents, so the bounds require additional work.
This fact manifests itself in our proof of \eqref{e.bd.Jop.e}.

Here is an overview of the proof.
As will be seen in Section~\ref{s.limiting}, proving Theorem~\ref{t.main} boils down to bounding the weighted norms of the various operators that make up $\sg(t)$ and $\sg^{\e}(t)$.
We obtain these bounds in two steps, by first employing a comparison argument to reduce proving these bounds to proving their unweighted analogs, and then proving their unweighted analogs.
While most of the unweighted analogs follow from existing results, one of them, \eqref{e.bd.Jop.e}, does not seem straightforward and requires significant work.
One may also consider the weighted $\Lsp^{p}$ space for $p\neq 2$, and we leave that to future work.

Let us mention a few more works on the moments of the noise-mollified SHE and point to the references in \cite{caravenna2023critical} for more references.
The work \cite{feng2015} carried out detailed analysis of the second moment and further used of the Gaussian correlation inequality to compare the moments of the SHE with those of a log-normal.
The works \cite{caravenna2023gmc,clark2023planar,clark2024continuum} studied properties related to the Gaussian multiplicative chaos.
High moments of the polymers and the mollified SHE are studied in \cite{cosco2023moments,lygkonis2023moments,cosco2024moments}.

The rest of the paper is organized as follows.
In Section~\ref{s.limiting}, we recall the definition of $\sg(t)$ and prove Theorem~\ref{t.main}\eqref{t.main.1}.
In Section~\ref{s.limiting}, we consider the prelimiting semigroup $\sg^{\e}(t):=e^{-t\ham^{\e}}$ and prove Theorem~\ref{t.main}\eqref{t.main.2}.

\section{The limiting semigroup, Proof of Theorem~\ref{t.main}\eqref{t.main.1}}
\label{s.limiting}
We begin by recalling the definition of the delta-Bose semigroup $\sg(t)$.
For $\alpha\in\pair\intv{n}$, the relevant operators map between functions on 
\begin{align}
	\label{e.xsp}
	\R^{2\intv{n}} 
	:=
	(\R^{2})^{\intv{n}}
	&:= 
	\big\{ (x_i)_{i\in\intv{n}} \, \big| \, x_i\in\R^2 \big\}
	\ ,
\\
	\label{e.ysp}
	\R^{2}\times\R^{2\intv{n}\setminus\alpha} 
	:= 
	\R^{2}\times\R^{2(\intv{n}\setminus\alpha)} 
	&:=
	\big\{ y=(\yc, (y_i)_{i\in\intv{n}\setminus\alpha}) \, \big| \, y_c,y_i\in\R^2 \big\}
	\ ,
\end{align}
where we index the first coordinate in \eqref{e.ysp} by ``c'' for ``center of mass''.
Consider
\begin{align}
	S_{\alpha}: \R^{2}\times\R^{2\intv{n}\setminus\alpha}  \to \R^{2\intv{n}},
	\qquad
	\big( S_{\alpha}y \big)_{i}
	:=
	\begin{cases}
		\yc & \text{when } i\in\alpha
	\\
		y_i & \text{when } i\in\intv{n}\setminus\alpha\ .
	\end{cases}
\end{align}
Let $\hk(t,x_i):=\exp(-\normtwo{x_i}^2/2t)/(2\pi t)$ denote the heat kernel on $\R^2$, let $\heatsg(t,x):=\prod_{i\in\intv{n}}\hk(t,x_i)$ denote the heat kernel on $\R^{2\intv{n}}$, and let
\begin{align}
	\label{e.jfn}
	\jfn(t)
	:=
	\int_0^{\infty} \d u \frac{t^{u-1}e^{\theta u}}{\Gamma(u)}\ ,
\end{align}
where $\theta\in\R$ is the parameter in \eqref{e.beta}.
For $\alpha\neq\alpha'\in\pair\intv{n}$, define the integral operators $\heatsg_{\alpha}(t)$, $\heatsg_{\alpha}(t)^*$, $\heatsg_{\alpha\alpha'}(t)$, $\Jop_{\alpha}(t)$ through their kernels as
\begin{subequations}
\label{e.ops}
\begin{align}
	\label{e.incoming}
	\heatsg_{\alpha}(t,y,x)
	&:=
	\heatsg\big(t, S_{\alpha}y - x \big)
	=:
	\big(\heatsg_{\alpha}\big)^*(t,x,y)
	\ ,
\\
	\label{e.swapping}
	\heatsg_{\alpha\alpha'}(t,y,y')
	&:=
	\heatsg\big(t, S_{\alpha}y - S_{\alpha'}y' \big)\ ,
\\
	\label{e.Jop}
	\Jop_{\alpha}(t,y, y')
	&:=
	4\pi\,\jfn(t) \, \hk(\tfrac{t}{2},\yc-\yc') \cdot \prod_{i\in\intv{n}\setminus\alpha} \hk(t,y_i-y'_i)\ ,
\end{align}
\end{subequations}
where $x\in\R^{2\intv{n}}$ and $y,y'\in\R^{2}\times\R^{2\intv{n}\setminus\alpha}$ in \eqref{e.incoming} and \eqref{e.Jop}, and $y\in\R^{2}\times\R^{2\intv{n}\setminus\alpha}$ and $y'\in\R^{2\intv{n}\setminus\alpha'}$ in \eqref{e.swapping}.
Next, let
\begin{align}
	\label{e.dgm}
	\dgm\intv{n}
	:=
	\big\{ \vecalpha=(\alpha_k)_{k=1}^m\in\pair\intv{n}^m \, \big| \, m\in\Z_{>0}, \alpha_{k}\neq\alpha_{k+1} \text{ for } k=1,\ldots, m-1 \big\}
	\ .
\end{align}
This set indexes certain diagrams, hence the name $\dgm$; see \cite[Section~2]{gu2021moments}.
Write $|\vecalpha|:=m$ for the length of $\vecalpha\in\dgm\intv{n}$.
For $f=f(\tau,\tau',\tau'',\ldots)$ that depends on finitely many nonnegative $\tau$s, write $\int_{\Sigma(t)}\d \vectau f = \int_{\tau+\tau'+\ldots=t} \d \vectau f$ for the convolution-like integral.
For $\vecalpha\in\dgm\intv{n}$, let
\begin{align}
	\label{e.sgsum}
	\sgsum_{\vecalpha}(t)
	&:=
	\int_{\Sigma(t)} \d \vectau \
	\heatsg_{\alpha_{1}}(\tau_{\frac{1}{2}})^*
	\prod_{k=1}^{|\vecalpha|-1} \Jop_{\alpha_{k}}(\tau_{k}) \, \heatsg_{\alpha_{k}\alpha_{k+1}}(\tau_{k+\frac{1}{2}}) \cdot
	\Jop_{\alpha_{|\vecalpha|}}(\tau_{|\vecalpha|}) \, \heatsg_{\alpha_{|\vecalpha|}}(\tau_{|\vecalpha|+\frac{1}{2}})\ .
\end{align}
Hereafter, products of operators are understood in the written order, so $\prod_{k=1}^{K}\Top_k := \Top_1 \Top_2 \cdots \Top_K$.
The delta-Bose semigroup on $\R^{2\intv{n}}$ is
\begin{align}
\label{e.sg}
	\sg(t)
	:=
	\heatsg(t)
	+
	\sum_{\vecalpha\in\dgm\intv{n}} 
	\sgsum_{\vecalpha}(t)\ .
\end{align}

In \cite[Section~8]{gu2021moments}, the series \eqref{e.sg} is shown to converge in the $\Lsp^2$ operator norm.
Here we seek to do similarly for the weighted operator norm.
Let us prepare some notation and tools.
First,
\begin{align}
	\label{e.norma.}
	\norma{f}^2
	&=
	\norm{f}_{\Lsp^{2}_{a}(\Omega)}^2
	:=
	\begin{cases}
		\int_{\Omega} \d x \, \big| f(x)e^{a\normone{x}} \big|^{2}
		&
		\Omega = \R^{2\intv{n}}\ ,
	\\
		\int_{\Omega} \d y\,  \big| f(y)e^{2a\normone{\yc}+a\sum_{i\in\intv{n}\setminus\alpha}\normone{y_i}} \big|^{2}
		&
		\Omega = \R^{2}\times\R^{2\intv{n}\setminus\alpha}\ .
	\end{cases}	
\end{align}
For $\Omega=\R^{2\intv{n}}$, the norm is the same as \eqref{e.norma}, but for $\Omega=\R^2\times\R^{2\intv{n}\setminus\alpha}$, we assign twice the exponential weight to $\yc$, which is natural because it bears a meaning of ``merging'' the two coordinates in $\alpha$.
Accordingly, for operators that map between functions on the spaces \eqref{e.xsp}--\eqref{e.ysp}, define the weighted operator norm $\normopa{\cdot}$ as \eqref{e.normopa} with the $\norma{\cdot}$ given in \eqref{e.norma.}.
We will omit the underlying spaces $\Omega$ when writing the operator norms, because the spaces can be read off from the definition of the operators.
For example, referring to \eqref{e.incoming} and the description that follows, we see that 
$
	\normopa{\heatsg_{\alpha}(t)}
$
means the operator norm of $\heatsg_{\alpha}(t):\Lsp^{2}_{a}(\R^{2\intv{n}})\to \Lsp^{2}_{a}(\R^2\times\R^{2\intv{n}\setminus\alpha})$.
Next, note that the kernels in \eqref{e.ops} are nonnegative.
For a one-parameter family of integral operators $\Top(t)$, $t>0$, with a nonnegative kernel $\Top(t,z,z')$, if
\begin{align}
	\sup\Big\{ \int_0^b \d t \int \d z' \d z \ |f'(z')| \cdot \Top(t,z',z) \cdot |f(z)| \Big\}
	\ ,
\end{align}
is finite when the supremum is taken over $\norma{f'}\leq1$ and $\norma{f}\leq 1$, then the operator $\int_0^b \d t\,\Top(t)$ is well-defined and bounded with $\Normopa{ \int_0^b \d t \, \Top(t) }$ equal to the above expression.
When $a=0$, the weighted norms reduce to the $\Lsp^2$ norms, which we write as $\norm{\cdot}_{2,0}=\norm{\cdot}_{2}$ and $\norm{\cdot}_{2,0\to 2,0}=\norm{\cdot}_{2\to 2}$. It will be useful to be able to bound convolution-like integral as in equation \eqref{e.sgsum} in norm by bounding the individual integrands. To this end, we have the following lemma:
\begin{lemma}[{\cite[Lemma~2.1(a)]{tsai2024}; also \cite[Lemma 8.10]{gu2021moments}}]
\label{l.sum}
For $m\in\Z_{>0}$ and $\kappa\in (\frac{1}{2}\Z)\cap(0,m+1)$, let $\Top_{\kappa}(t):\Bsp_{\kappa} \to \Bsp_{\kappa-1/2}$ be a bounded operator with a nonnegative kernel, where $\Bsp_{\kappa}$ is a Banach space consisting of some Borel functions on $\R^{d_\kappa}$, and let $\norm{\Top_{\kappa}(t)}_{\mathrm{op}}$ denote the operator norm.
Assume that, for all $t>0$ and a constant $c_0\in(0,\infty)$,
\begin{align}
	\label{e.l.sum.bd}
	\Norm{ \Top_{\kappa}(t) }_{\mathrm{op}}
	&\leq
	c_0 e^{c_0t}
	\cdot
	\begin{cases}
		t^{-1/2} & \text{when } \kappa = \tfrac{1}{2} \ , 
	\\
		t^{-1}\, \big|\log(\tfrac{1}{2}\wedge t)\big|^{-2} & \text{when } \kappa \in \Z\cap [1,m] \ , 
	\\
		t^{-1} & \text{when } \kappa \in (\tfrac{1}{2}+\Z)\cap(1,m) \ , 
	\\
		t^{-1/2} & \text{when } \kappa = m+ \tfrac{1}{2} \ , 
	\end{cases}
\\
	\label{e.l.sum.bdint}
	\NOrm{ \int_0^{\infty} \d t \, e^{-c_0 t} \Top_{\kappa}(t) }_{\mathrm{op}}
	&\leq
	c_0\ 
	\text{ when } \kappa \in (\tfrac{1}{2}+\Z)\cap(1,m) \ .
\end{align}
Fix $c_1 \geq 0$ and let $\Top'_{\kappa}(t):=c_1\delta_{0}(t)\ind+\Top_{\kappa}(t)$ when $\kappa\in\Z\cap[1,m]$ and $\Top'_{\kappa}(t):=\Top_{\kappa}(t)$ when $\kappa\in(\frac12+\Z)\cap(0,m+1)$.
Then, there exists a universal $c\in(0,\infty)$, such that for all $m\in\Z_{>0}$, $t>0$, and $\lambda \geq c_0+2$,
\begin{align}
	\NOrm{ \int_{\Sigma(t)} \d \vectau \, \prod_{k=1}^{2m+1} \Top'_{k/2}(\tau_{k/2}) }_{\mathrm{op}}
	\leq
	c^m \, m^3 \, e^{\lambda t}
	\Big( c_1^m+c_0^{2}\Big(c_1+\frac{c_0^2}{\log(\lambda-c_0-1)}\Big)^{m-1} \Big)\ .
\end{align}
\end{lemma}

Given Lemma~\ref{l.sum}, our task is reduced to bounding the norm $\normopa{\cdot}$ of the operators in \eqref{e.ops}. The Lemma then implies the convergence of the operator in equation \eqref{e.sg} by choosing $\lambda$ large enough. First, consider the multiplicative isometric operator $\weiop:\Lsp^{2}_{a}(\R^{2})\to \Lsp^{2}(\R^{2})$, $\psi\mapsto e^{a\normone{\cdot}}\psi$ and express the weights in the norm as conjugations:
\begin{subequations}
\label{e.pa.p}
\begin{align}
	\label{e.pa.p.incoming}
	\Normopa{ \heatsg_{\alpha}(t) }
	&=
	\Normop{ \weiop^{2}\otimes\weiop^{\otimes\intv{n}\setminus\alpha} \cdot \heatsg_{\alpha}(t) \cdot \weiop^{-\otimes\intv{n}} }\ ,
\\
	\Normopa{ \heatsg_{\alpha}(t)^* }
	&=
	\Normop{ \weiop^{\otimes\intv{n}} \cdot \heatsg_{\alpha}(t)^* \cdot \weiop^{-2}\otimes\weiop^{-\otimes\intv{n}\setminus\alpha} }\ ,
\\
	\Normopa{ \heatsg_{\alpha\alpha'}(t) }
	&=
	\Normop{ 
		\weiop^{2}\otimes\weiop^{\otimes\intv{n}\setminus\alpha}
		\cdot 
		\heatsg_{\alpha\alpha'}(t) 
		\cdot 
		\weiop^{-2}\otimes\weiop^{-\otimes\intv{n}\setminus\alpha'} 
	}\ ,
\\
	\NOrmopa{ \int_0^\infty \d t \, e^{-\gamma t} \heatsg_{\alpha\alpha'}(t) }
	&=
	\NOrmop{ 
		\weiop^{2}\otimes\weiop^{\otimes\intv{n}\setminus\alpha}
		\cdot
		\int_0^\infty \d t \, e^{-\gamma t} \heatsg_{\alpha\alpha'}(t)
		\cdot
		\weiop^{-2}\otimes\weiop^{-\otimes\intv{n}\setminus\alpha'} 
	}\ ,
\\
	\Normopa{ \Jop_{\alpha}(t) }
	&=
	\Normop{ \weiop^{2}\otimes\weiop^{\otimes\intv{n}\setminus\alpha} \cdot \Jop_{\alpha}(t) \cdot \weiop^{-2}\otimes\weiop^{-\otimes\intv{n}\setminus\alpha} }\ .
\end{align}
\end{subequations}
Next, note that the heat kernel satisfies the bound
\begin{align}
	\label{e.hk.comp}
	\hk(t,x_i-y_i) e^{a \,(\normone{y_i}-\normone{x_i})}
	\leq
	2 e^{2at}\hk(2t,x_i-y_i)\ ,
	\qquad
	x_i, y_i\in\R^2\ .
\end{align}
To see how, use the triangle inequality followed by Young's inequality to get
$
	\normone{y_i}-\normone{x_i} 
	\leq 
	\normone{x_i-y_i}
	\leq
	2a t + \normone{x-y}^2/(8a t) 
	\leq 
	2a t +\normtwo{x-y}^2/(4a t)
$,
where we used the bound $\normone{\cdot}^2 \leq 2 \normtwo{\cdot}^2$ that holds on $\R^{2}$.
This bounds the left-hand side of \eqref{e.hk.comp} by
$
	\hk(t,x_i-y_i) \, e^{ 2at + \normtwo{x_i-y_i}^2/(4t) }
$.
Further using the readily verified identity $\hk(t,x_i-y_i)e^{\normtwo{x_i-y_i}^2/(4t)} = 2 \hk(2t,x_i-y_i) $ gives \eqref{e.hk.comp}.
Referring to \eqref{e.incoming}, we see that the conjugated operator in \eqref{e.pa.p.incoming} has kernel $e^{a|S_\alpha y|_1-a|x|_1}\heatsg(t,S_\alpha y-x)$.
Using \eqref{e.hk.comp} and the fact that the kernel is nonnegative, we bound the right-hand side of \eqref{e.pa.p.incoming} by $ce^{ct}\normop{\heatsg_{\alpha}(2t)}$.
Hereafter, we write $c=c(\theta,n,a,\Phi)$ for a general, finite, positive, deterministic constant that may change from place to place but depends only on $\theta,n,a,\Phi$.
Similar arguments apply to the operators in \eqref{e.pa.p} and give, for $\alpha\neq\alpha'$, the comparison bounds
\begin{subequations}
\label{e.comp}
\begin{align}
	\label{e.comp.incoming}
	\Normopa{ \heatsg_{\alpha}(t) }
	&\leq
	ce^{ct}
	\Normop{ \heatsg_{\alpha}(2t)  }\ ,
\\
	\label{e.comp.incoming*}
	\Normopa{ \heatsg_{\alpha}(t)^* }
	&\leq
	ce^{ct}
	\Normop{ \heatsg_{\alpha}(2t)^*  }\ ,
\\
	\label{e.comp.swapping}
	\Normopa{ \heatsg_{\alpha\alpha'}(t) }
	&\leq
	ce^{ct}
	\Normop{ \heatsg_{\alpha\alpha'}(2t) }\ ,
\\
	\label{e.comp.swapping.int}
	\NOrmopa{ \int_0^\infty \d t \, e^{-(c+2)t} \heatsg_{\alpha\alpha'}(t) }
	&\leq
	c\,
	\NOrmop{ \int_0^\infty \d t \, e^{-2t}\heatsg_{\alpha\alpha'}(2t) }\ ,
\\
	\label{e.comp.Jop}
	\Normopa{ \Jop_{\alpha}(t) }
	&\leq
	ce^{ct} \jfn(t)\, \normop{\hk(t)}\cdot \normop{\hk(2t)}^{n-2}\ .
\end{align}
\end{subequations}

We are now ready to prove Theorem~\ref{t.main}\eqref{t.main.1}.
By \cite[Lemmas~ 5.1, 8.1(a), 8.2, and 8.4]{gu2021moments} and the property that $\normop{\hk(s)}\leq 1$, the left-hand sides of \eqref{e.comp} are bounded by
$ c e^{ct}t^{-1/2} $,
$ c e^{ct}t^{-1/2} $,
$ c e^{ct}t^{-1} $,
$ c e^{ct} $, and
$ c e^{ct}t^{-1} |\log(\frac{1}{2}\wedge t)|^{-2} $, respectively.
Combining these bounds with Lemma~\ref{l.sum} for a large enough $\lambda$ gives Theorem~\ref{t.main}\eqref{t.main.1}.

\section{The prelimiting semigroup, proof of Theorem~\ref{t.main}\eqref{t.main.2}}
\label{s.prelimiting}
Let us introduce the $\e$ analogs of the operators in Section~\ref{s.limiting}.
For $\alpha=ij\in\pair\intv{n}$, 
\begin{align}
	\label{e.ysp.}
	&\R^{4}\times(\R^2)^{\intv{n}\setminus\alpha} 
	:= 
	\R^{4}\times\R^{2\intv{n}\setminus\alpha} 
	:= 
	\big\{ y=(\yr, \yc, (y_i)_{i\in\intv{n}\setminus\alpha}) \, \big| \, \yr, \yc,y_i\in\R^2 \big\}
	\ ,
\\
	&S^\e_{\alpha}: \R^{4}\times\R^{2\intv{n}\setminus\alpha} \to \R^{2\intv{n}},
	\qquad
	\big( S^\e_{\alpha} y \big)_{k}
	:=
	\begin{cases}
		\yc+\e\yr/2 & \text{when } k=i
	\\
		\yc-\e\yr/2 & \text{when } k=j
	\\
		y_{k} & \text{when } k\in\intv{n}\setminus\alpha\ .
	\end{cases}
\end{align}
where we index the first two coordinates in \eqref{e.ysp.} respectively by $\textr$ and $\textc$ for ``relative'' and  ``center of mass''.
Set $\phi:=\sqrt{\Phi}$ and view $\Phi$ and $\phi$ as multiplicative operators acting on $\Lsp^2(\R^{2})$, where $\R^{2}$ represents the space of the $\yr$ coordinate.
For $\alpha,\alpha'\in\pair\intv{n}$, define operators 
\begin{subequations}
\label{e.ops.e}
\begin{align}
	\label{e.incoming.e}
	\heatsg^{\e}_{\alpha}(t,y,x)
	&:=
	\phi(\yr)\,
	\heatsg\big(t, S^{\e}_{\alpha}y - x \big)
	=:
	\big(\heatsg^{\e}_{\alpha}\big)^*(t,x,y)\ ,
\\
	\label{e.swapping.e}
	\heatsg^{\e}_{\alpha\alpha'}(t,y,y')
	&:=
	\phi(\yr)\,
	\heatsg \big(t, S^{\e}_{\alpha}y - S^{\e}_{\alpha'}y' \big)\,
	\phi(\yr')\ ,
\\	
	\label{e.Jop.e}
	\Jop^{\e}_{\alpha}(t)
	&:=
	\sum_{k=1}^\infty \beta_\e^{k+1}
	\int_{\Sigma(t)} \d\vectau\,
	\heatsg^{\e}_{\alpha\alpha}(\tau_1)\cdots\heatsg^{\e}_{\alpha\alpha}(\tau_k)\ ,
\\
\begin{split}
	\sgsum^{\e}_{\vecalpha}(t)
	&:=
	\int_{\Sigma(t)} 
	\d \vectau \,
	\heatsg^{\e}_{\alpha_{1}}(\tau_{\frac{1}{2}})^*
	\prod_{k=1}^{|\vecalpha|-1} 
	\Big( \beta_\e \delta_{0}(\tau_k) + \Jop^{\e}_{\alpha_{k}}(\tau_{k}) \Big) \, 
	\heatsg^{\e}_{\alpha_{k}\alpha_{k+1}}(\tau_{k+\frac{1}{2}})
\\
	&\hspace{.1\linewidth}\cdot
	\Big( \beta_\e \delta_{0}(\tau_{|\vecalpha|}) + \Jop^{\e}_{\alpha_{|\vecalpha|}}(\tau_{|\vecalpha|}) \Big)
	\, 
	\heatsg^{\e}_{\alpha_{|\vecalpha|}}(\tau_{|\vecalpha|+\frac{1}{2}})
	\ ,
\end{split}
\end{align}
\end{subequations}
where $x\in\R^{2\intv{n}}$, $y\in\R^{4}\times\R^{2\intv{n}\setminus\alpha}$, and $y'\in\R^{4}\times\R^{2\intv{n}\setminus\alpha'}$.
The operator $\Jop^{\e}_{\alpha}(t)$ permits another expression.
Let $\hk_\e(t,\yr-\yr'):=\hk(t,\e(\yr-\yr'))$ and consider
\begin{align}
	\label{e.Jop.e.}
	\jfn^{\e}(t):\Lsp^2(\R^2)\to\Lsp^2(\R^2)\ ,
	\quad
	\jfn^{\e}(t)
	:=
	\sum_{\ell=1}^\infty \beta_\e^{\ell+1}
	\int_{\Sigma(t)} \d\vectau\,
	\prod_{k=1}^{\ell} \phi\,\hk_{\e}(2\tau_{k})\,\phi.
\end{align}
It is straightforward to check from \eqref{e.Jop.e} that
\begin{align}
	\Jop^{\e}_{\alpha}(t,y,y')
	=
	\jfn^{\e}(t,\yr,\yr') \, \hk(\tfrac{t}{2},\yc-\yc') \, \prod_{i\in\intv{n}\setminus\alpha} \hk(t,y_i-y'_i)\ .
\end{align}

The prelimiting semigroup $\sg^{\e}(t):=e^{-t\ham^{\e}}$ enjoys an expansion similar to that of $\sg(t)$:
\begin{align}
\label{e.sg.e}
	\sg^{\e}(t)
	=
	\heatsg(t)
	+
	\sum_{\vecalpha\in\dgm\intv{n}} 
	\sgsum^{\e}_{\vecalpha}(t)\ .
\end{align}
To prove this, in \eqref{e.ham}, use Duhamel's principle to write $\sg^{\e}(t)=\heatsg(t)+\int_0^t \d s\,\heatsg(t-s)\sum_{\alpha}\Phi^{\e}_{\alpha}\sg^{\e}(t)$ and iterate the equation to obtain
\begin{align}
	\label{e.duhamel}
	\sg^{\e}(t)
	=
	\heatsg(t)
	+
	\sum_{\veceta} \int_{\Sigma(t)} \d \tau\,
	\heatsg^{\e}_{\eta_1}(\tau_0)^* \prod_{k=2}^{|\veceta|}\heatsg^{\e}_{\eta_{k-1}\eta_{k}}(\tau_{k}) \cdot \heatsg^{\e}_{\eta_{|\veceta|}}(\tau_{|\veceta|})\ ,
\end{align}
where the sum runs over $\veceta\in\cup_{m=1}^\infty(\pair\intv{n})^m$.
Rewrite $\veceta$ as $(\alpha_1^{k_1},\alpha_2^{k_2},\ldots)$, where $\alpha^k:=(\alpha,\ldots,\alpha)\in(\pair\intv{n})^k$, and $\alpha_1\neq\alpha_2$, $\alpha_2\neq\alpha_3$, \ldots.
Accordingly, the sum in \eqref{e.duhamel} is rewritten as sums over $\vecalpha\in\dgm\intv{n}$ and over $k_1,\ldots,k_{|\vecalpha|}\in\Z_{>0}$.
Carrying out the latter sum gives \eqref{e.sg.e}.

Next, similarly define the weighted norm $\norma{f}=\norm{f}_{\Lsp^{2}_{a}(\Omega)}$ by
\begin{align}
	\label{e.norma.e}
	\norm{f}_{\Lsp^{2}_{a}(\Omega)}^2
	:=
	\begin{cases}
		\int_{\Omega} \d x \, \big| f(x)e^{a|x|_1} \big|^{2}
		& 
		\Omega = \R^{2\intv{n}}\ ,
	\\
		\int_{\Omega} \d y \, \big| f(y)e^{ a\sum_{\sigma=\pm 1} \normone{\yc+\frac{\e\sigma}{2}\yr} + a \sum_{i\in\intv{n}\setminus\alpha} \normone{y_i} } \big|^{2}
		& 
		\Omega = \R^{4}\times\R^{2\intv{n}\setminus\alpha}\ .
	\end{cases}	
\end{align}
We have slightly abused notation by using the same notation for the norms in \eqref{e.norma} and in \eqref{e.norma.e}.
This should not cause any confusion, since the latter will only be applied to $\e$-dependent operators.
Define the operator norm $\normopa{\cdot}$ the same way as in \eqref{e.normopa} with the norm in \eqref{e.norma.e} replacing that of \eqref{e.norma}.

The comparison argument in Section~\ref{s.limiting} applies here, too.
Consider the multiplicative operator $\til{\weiop}:\Lsp^{2}_{a}(\R^{4})\to \Lsp^{2}_{a}(\R^{4})$, $(\til{\weiop}\psi)(\yr,\yc):= e^{a\sum_{\sigma=\pm 1} \normone{\yc+\frac{\e\sigma}{2}\yr} }\psi(\yr,\yc)$.
The bounds in \eqref{e.pa.p} continue to hold when we replace $\heatsg_{\alpha}(t)$, $\heatsg_{\alpha}(t)^*$, $\heatsg_{\alpha\alpha'}(t)$, $\Jop_{\alpha}(t)$, $\weiop^2$, and $\weiop^{-2}$ with $\heatsg^{\e}_{\alpha}(t)$, $\heatsg^{\e}_{\alpha}(t)^*$, $\heatsg^{\e}_{\alpha\alpha'}(t)$, $\Jop^{\e}_{\alpha}(t)$, $\til{\weiop}$, and $\til{\weiop}^{-1}$, respectively.
Next, the same comparison argument leading to \eqref{e.comp} applies also to their $\e$ analogs and gives, for $\alpha\neq\alpha'$,
\begin{subequations}
\label{e.comp.e}
\begin{align}
	\Normopa{ \heatsg^{\e}_{\alpha}(t) }
	&\leq
	ce^{ct}
	\Normop{ \heatsg^{\e}_{\alpha}(2t)  }\ ,
\\
	\Normopa{ \heatsg^{\e}_{\alpha}(t)^* }
	&\leq
	ce^{ct}
	\Normop{ \heatsg^{\e}_{\alpha}(2t)^*  }\ ,
\\
	\Normopa{ \heatsg^{\e}_{\alpha\alpha'}(t) }
	&\leq
	ce^{ct}
	\Normop{ \heatsg^{\e}_{\alpha\alpha'}(2t)  }\ ,
\\
	\NOrmopa{ \int_0^\infty \d t \,e^{-(c+2)t} \heatsg^{\e}_{\alpha\alpha'}(t) }
	&\leq
	c\,
	\NOrmop{ \int_0^\infty \d t \, e^{-2t}\heatsg^{\e}_{\alpha\alpha'}(2t) }\ .
\\
	\label{e.comp.Jop.e}
	\Normopa{ \Jop^{\e}_{\alpha}(t) }
	&\leq
	ce^{ct} \, \normopa{\jfn^{\e}(t)} \, \normop{\hk(t)}\cdot \normop{\hk(2t)}^{n-2}\ .
\end{align}
\end{subequations}
Further, since $\phi$ has a compact support, $\normopa{\jfn^{\e}(t)}\leq c\normop{\jfn^{\e}(t)}$.
Using this bound and $\normop{\hk(s)}\leq 1$ turns \eqref{e.comp.Jop.e} into
$
	\normopa{ \Jop^{\e}_{\alpha}(t) }
	\leq
	ce^{ct} \, \normopa{\jfn^{\e}(t)}
$.
Given these comparison bounds, the main task is to prove the following bounds.
\begin{subequations}
\label{e.bds.e}
\begin{align}
	\label{e.bd.incoming.e}
	\Normop{ \heatsg^{\e}_{\alpha}(t) }
	=
	\Normop{ \heatsg^{\e}_{\alpha}(t)^* }
	&\leq
	c \, t^{-1/2}\ ,
\\
	\label{e.bd.swapping.e}
	\Normop{ \heatsg^{\e}_{\alpha\alpha'}(t) }
	&\leq
	c \, t^{-1}\ ,
\\
	\label{e.bd.swapping.int.e}
	\NOrmop{  \int_0^\infty \d t\, e^{-t}\heatsg^{\e}_{\alpha\alpha'}(t) }
	&\leq
	c\ ,
\\
	\label{e.bd.Jop.e}
	\Normop{ \jfn^{\e}(t) }
	&\leq
	c \, t^{-1} \big|\log\big(t\wedge \tfrac{1}{2}\big)\big|^{-2} \cdot e^{c\,t}\ .
\end{align}
\end{subequations}

Postponing the proof of \eqref{e.bds.e} to Section~\ref{s.prelimiting.bds}, we now complete the proof of Theorem~\ref{t.main}\eqref{t.main.2}.
First, combining \eqref{e.bds.e} and Lemma~\ref{l.sum} for $c_1=\beta_\e$ and for a large enough $\lambda$ gives the first statement in Theorem~\ref{t.main}\eqref{t.main.2}.
To prove the second statement, let $(\ind_{\leq M}f)(x):=f(x)\ind_{\normone{x}\leq M}$, let $\ind_{>M}:=\ind-\ind_{\leq M}$, and write
\begin{align}
\begin{split}
	\Norma{ &(\sg(t)-\sg^{\e}(t))f }
	\leq
	\Norma{ \ind_{\leq M}(\sg(t)-\sg^{\e}(t))\ind_{\leq M'} f }
\\
	&+
	\sum_{\Top=\sg(t),\sg^{\e}(t)}
	\Norma{ \ind_{> M}\Top \ind_{\leq M'} f }
	+
	\sum_{\Top=\sg(t),\sg^{\e}(t)}
	\Norma{ \Top \ind_{>M'} f }\ . 
\end{split}
\end{align}
On the right-hand side, send $\e \to 0$ first, $M\to\infty$ second, and $M'\to\infty$ last.
The first term on the right-hand side is bounded by $e^{|a|M}\norm{(\sg(t)-\sg^{\e}(t))\ind_{\leq M'} f }_2$.
For fixed $M'$, $\ind_{\leq M'} f \in\Lsp^2(\R^{2\intv{n}})$, so by \cite[Theorem~1.6(b)]{gu2021moments} the term tends to $0$ as $\e\to 0$.
The second term on the right-hand side is bounded by $e^{-M}\norm{\sg(t) \ind_{\leq M'}f}_{2,a+1}+e^{-M}\norm{\sg^{\e}(t) \ind_{\leq M'}f}_{2,a+1}$, which is in turn bounded by $e^{-M} \, (\norm{\sg(t)}_{2,a+1\to 2,a+1}+\norm{\sg^{\e}(t)}_{2,a+1\to 2,a+1})\, \norm{\ind_{\leq M'}f}_{2,a+1}$.
This tends to $0$ as $\e\to 0$ first and $M\to\infty$ second.
The last term on the right-hand side is bounded by $(\norm{\sg(t)}_{2,a\to 2,a}+\norm{\sg^{\e}(t)}_{2,a\to 2,a})\, \norm{\ind_{>M'}f}_{2,a}$, which tends to $0$ as $\e\to 0$ first and $M'\to\infty$ later.

\subsection{Proving~\eqref{e.bds.e}}
\label{s.prelimiting.bds}
The bounds \eqref{e.bd.incoming.e}--\eqref{e.bd.swapping.int.e} follow from existing results.
For $\alpha=ij$, let $N_{\alpha}$ denote the map $\R^{2\intv{n}}\to \R^{2}\times\R^{2\intv{n}\setminus\alpha}, N_{\alpha}p:=(p_i+p_j, p_{\intv{n}\setminus\alpha})$, and let $p^{-}_{\alpha}:=(p_i-p_j)/2$.
For $f\in\Lsp^2(\R^{2\intv{n}})$ and $g\in\Lsp^2(\R^{4}\times\R^{2\intv{n}\setminus\alpha})$, let
\begin{align}
	\Fou{f}(p)
	:=
	\int_{\R^{2\intv{n}}} \frac{\d p}{(2\pi)^{n}} e^{-\img p\cdot x} f(x)\ , 
	&&
	\pFou{g}(\yr, q)
	:=
	\prod_{i\in \textc\cup \intv{n}\setminus\alpha}
	\int_{\R^2} \frac{\d y_i}{2\pi} e^{-\img q_i\cdot y_i}
	\cdot
	g(y)\ 
\end{align}
denote the Fourier and partial Fourier transforms, where $p\in\R^{2\intv{n}}$ and $q\in \R^{2}\times\R^{2\intv{n}\setminus\alpha}$.
By \cite[Equations~(4-11), (5-5b')]{gu2021moments},
\begin{align}
	\label{e.fourier.incoming}
	\Ip{ g, \heatsg^{\e}_{\alpha}(t) f }
	&=
	\int_{\R^2} \frac{\d \yr}{2\pi}	
	\int_{\R^{2\intv{n}}} \d p
	\,
	\bar{ \pFou{g}(\yr,N_{\alpha}p) } 
	\,
	\phi(\yr) e^{\img \e\yr\cdot p^{-}_{\alpha}}
	\,
	e^{-\frac{t}{2}\normtwo{p}^2}
	\Fou{f}(p)\ ,
\\
	\label{e.fourier.swapping}
\begin{split}	
	\Ip{ g_1, \heatsg^{\e}_{\alpha\alpha'}(t) g_2 }
	&=
	\int_{\R^2} \frac{\d \yr}{2\pi}
	\int_{\R^2} \frac{\d \yr'}{2\pi}
	\int_{\R^{2\intv{n}}} \d p
\\
	&
	\cdot
	\bar{ \pFou{g_1}(\yr,N_{\alpha} p) } 
	\,
	\phi(\yr)
	\,
	e^{-\frac{t}{2}\normtwo{p}^2+\img \e(\yr\cdot p^{-}_\alpha-\yr'\cdot p^{-}_{\alpha'})} \phi(\yr')
	\,
	\pFou{g_2}(\yr',N_{\alpha'}p)\ .
\end{split}
\end{align}
In \eqref{e.fourier.incoming}, bound $|\phi(\yr)e^{\img\e\yr\cdot p^-_{\alpha}}|\leq c\,\phi(\yr)$ and $e^{-t\normtwo{p}^2/2}\leq e^{-t\normtwo{p_i-p_j}^2/4}$, write the resulting integrand as the product of $c\,\phi(\yr)\,\bar{ \pFou{g}(\yr,N_{\alpha}p) }e^{-t\normtwo{p_i-p_j}^2/4}$ and $\Fou{f}(p)$, and apply the Cauchy--Schwarz inequality over the $p$ integral.
Doing so gives the bound $c\,t^{-1/2}\,(\int_{\R^2} \d \yr \phi(\yr) \norm{g(\yr,\cdot)}_2)\cdot \norm{f}_2$.
Further applying the Cauchy--Schwarz inequality over the $\yr$ integral gives \eqref{e.bd.incoming.e}.
Applying a similar argument to \eqref{e.fourier.swapping} gives \eqref{e.bd.swapping.e}.
Applying $\int_0^\infty \d t\, e^{-t}$ to \eqref{e.fourier.swapping} followed by using \cite[Equations (3.1), (3.3), (3.4), (3.6)]{dellantonio1994hamiltonians} gives \eqref{e.bd.swapping.int.e}.

Proving~\eqref{e.bd.Jop.e} requires more work.
We begin with a preliminary bound.
Let $\norm{\cdot}_{\mathrm{HS}}$ denote the Hilbert--Schmidt norm and let $b(u):= \int_{\R^4} \d\yr\d\yr'\,\phi^2(\yr)e^{-u\normtwo{\yr-\yr'}^2/4}\phi^2(\yr')$. 
We have $\normop{\phi\hk_\e(2\tau)\phi} \leq \norm{\phi\hk_\e(2\tau)\phi}_{\mathrm{HS}} = b(\e^2/\tau)^{1/2}/4\pi\tau$.
Recall that $\phi$ is $\Csp^\infty$ smooth, has a compact support, and satisfies $\int_{\R^2}\d x\,\phi^2(x)=1$.
Using these properties, it is not hard to check that $b$ is $\Csp^\infty$ smooth and strictly decreasing on $u\in[0,\infty)$, that $b(0)=1$, that $b'(0)<0$, and that $\limsup_{u\to\infty}b(u)u < \infty$.
These properties together imply that $b(u) \leq (1+c_0u)^{-1}$, for some $c_0\in(0,\infty)$ depending only on $\phi$.
Hence $\normop{\phi\hk_\e(2\tau)\phi} \leq (4\pi)^{-1} (\tau^2+c_0\e^{2}\tau)^{-1/2}$.
Using this in \eqref{e.Jop.e.} gives
\begin{align}
	\label{e.bdingJop.B}
	\normopa{\jfn^{\e}(t)}
	\leq
	B_\e(t)\ ,
	\qquad
	B_\e(t) := \beta_\e \sum_{\ell=1}^\infty \int_{\Sigma(t)}\d\vectau\, \prod_{k=1}^{\ell} \frac{\beta_\e}{4\pi\sqrt{\tau_k^2+c_0\e^2\tau_k}}\ .
\end{align}

The bound \eqref{e.bdingJop.B} reduces our task to bounding $B_\e(t)$.
Let us consider the cases $t\leq \e^{2}c_0$ and $t>\e^{2}c_0$ separately.
For $t\leq \e^{2}c_0$, bound $1/\sqrt{\tau_k^2+c_0\e^2\tau_k} \leq 1/\sqrt{\e^2\tau_k}$ and use the Dirichlet integral formula
$
	\int_{\Sigma(t)} \d \vectau \, \prod_{k=1}^{\ell} \tau_{k}^{-1/2} = t^{\ell/2-1} \Gamma(1/2)^{\ell}/\Gamma(\ell/2)
$
to get
\begin{align}
	B_\e(t)
	\leq
	c\,\beta_\e\, t^{-1}\, \, \sum_{\ell=1}^\infty \frac{(c\beta_\e \sqrt{t\e^{-2}} )^{\ell}}{\Gamma(\ell/2)}
	=
	c\,\beta_\e\, t^{-1/2}\e^{-1}\, \, \sum_{\ell=1}^\infty \frac{(c\beta_\e)^{\ell} (\sqrt{t\e^{-2}})^{\ell-1}}{\Gamma(\ell/2)}\ .	
\end{align}
The series converges for all small enough $\e$ because $\beta_\e\to 0$.
Using $t\e^{-2} \leq c_0$ and $\beta_\e \leq c/|\log\e|$, we bound $B_\e(t)$ by $c\,t^{-1/2} \cdot \e^{-1}|\log\e|^{-2}$.
Since $t^{-1/2}(\log t)^{-2}$ is decreasing for $t \leq \e^{2}c_0$ as $\e \to 0$, the last expression is bounded by $c\,t^{-1/2} \cdot t^{-1/2} |\log t|^{-2}$.

Moving on to the case $t > \e^{2}c_0$, we begin by deriving a contour-integral formula of $B_\e(t)$.
Consider the Laplace transform of the factor in the product in \eqref{e.bdingJop.B}:
\begin{align}
	\label{e.laplace.C}
	C_\e(\lambda)
	:=
	\int_0^\infty \d s \, \frac{e^{-\lambda s} \beta_\e}{4\pi\sqrt{s^2+c_0\e^2s}}\ .
\end{align}
Perform a change of variables $\lambda s\mapsto s$, use integration by parts\\{}$\int_0^\infty\d s \, e^{-s}/\sqrt{s^2+ss_0}=-2\log\sqrt{s_0}+2\int_0^\infty \d s \, e^{-s}\log(\sqrt{s}+\sqrt{s+s_0})$, and use \eqref{e.beta} to simplify the result.
Doing so gives
\begin{align}
	\label{e.laplace.C.}
	C_\e(\lambda) 
	= 
	1 - \frac{\beta_\e}{4\pi}( \log\lambda - \eta_\e - D(c_0\e^2\lambda) )\ ,
\end{align}
where $\eta_\e$ is a real constant such that $\eta_\e\to \eta \in\R$ as $\e\to 0$, and
\begin{align}
	\label{e.laplace.D}
	D(\lambda) 
	:=
	2\int_0^\infty \d s \, e^{-s}\log(\sqrt{s}+\sqrt{s+\lambda})\ .
\end{align}
Put $c_1:=e^{\gamma+D(0)+1}$.
By \eqref{e.laplace.C.}, for all real $\lambda>c_1$ and small enough $\e$, $C_\e(\lambda)<1$.
By \eqref{e.laplace.C}, for all real $\lambda>0$, $C_\e(\lambda)>0$.
Since the integral in \eqref{e.bdingJop.B} is a sum of convolutions, the Laplace transform of $B_\e(t)$ is a geometric series 
$
	\int_0^\infty \d t \, e^{-t\lambda} B_\e(t)
	=
	\beta_\e \sum_{k=1}^\infty C_\e(\lambda)^k
	=
	\beta_\e C_\e(\lambda)/(1-C_\e(\lambda)),
$
for all $\lambda>c_1$.
The last expression analytically continues to $z\in \C\setminus (-\infty,c_1]$.
Given the Laplace transform, Mellin's inversion formula yields
$
	B_\e(t)
	=
	\int_{\supset} \frac{\d z}{2\pi\img} \, e^{zt} \beta_\e C_\e(z)/(1-C_\e(z))
$,
where $\supset\,:=\{(x-\img 0)| -\infty < x \leq c_1\}\cup\{(x+\img 0)| -\infty < x \leq c_1\}$ denotes the contour that goes from $-\infty$ to $c_1$ along the lower side of the real axis and then from $c_1$ to $-\infty$ along the upper side of the real axis.
Add $\int_{\supset}\frac{\d z}{2\pi\img} \, e^{zt} \beta_\e $, which is zero, to the right-hand side.
Doing so gives the contour-integral formula
\begin{align}
	\label{e.laplace.B}
	B_\e(t)
	=
	\int_{\supset} \frac{\d z}{2\pi\img} \, \frac{ e^{zt} \beta_\e }{1-C_\e(z)}	
	=
	\int_{\supset} \frac{\d z}{2\pi\img} \, \frac{e^{zt}\, 4\pi}{-\log z + \eta_\e + D(c_0\e^2z)}\ .
\end{align}

Let us bound \eqref{e.laplace.B}.
Express the last integral in \eqref{e.laplace.B} in real variables as
$
	B_\e(t)
	=
	\int_{-c_1}^{\infty} \d x \, 
	e^{-tx}\, 4\, G_\e(x) / (F_\e(x)^2+ G_\e(x)^2)
$,
where
\begin{align}
	F_\e(x) &:= \log |x| - \eta_\e - \Re(D(-c_0\e^2x+\img 0))\ ,
\\
	G_\e(x) &:= \big( \pi - \Im(D(-c_0\e^2x+\img 0)) \big)\ind_{x>0}\ .
\end{align}
Let $\log_+ u := \log(u\vee 1)$ and $u_+:=u\vee 0$ for $u\in\R$.
It is not hard to check from \eqref{e.laplace.D} that, for all $x\in[-c_1,\infty)$, $|G_\e(x)| \leq c$ and $F_\e(x) \geq \log |x| - c - \log_+( c_0\e^2 x)$.
Hence
\begin{align}
	B_\e(t)
	\leq
	c
	\int_{-c_1}^{\infty} \d x \, \frac{ e^{-tx} }{ (\log(|x|/c_2) - \log_+( c_0\e^2 x))_+{}^2 }\ ,
\end{align}
for some constant $c_2\in[1,\infty)$.
Divide the integrals into
\begin{align}
	\int_{-c_1}^{c_2e^3} + \int_{c_2e^3}^{(c_2e^3)\vee(1/t)} + \int_{(c_2e^3)\vee(1/t)}^{1/c_0\e^2} + \int_{1/c_0\e^2}^{\infty}
\end{align}
and call the results $I_1$ through $I_4$, respectively.
The integral $I_1$ is bounded by $ce^{c_1t}$.
For $I_2$, since $x \leq 1/t \leq 1/\e^{2}c_0$, the denominator of the integrand is $|\log(x/c_2)|^2$.
Forgoing the exponential gives $I_2\leq I_2':= \int_{c_2e^3}^{(c_2e^3)\vee (1/t)} \d x /(\log (x/c_2))^2$.
Integrate by parts to get
\begin{align}
	I'_2 
	= 
	\frac{(x-c_2e^{3})}{(\log x/c_2)^2}\Big|^{(c_2e^3)\vee (1/t)}_{c_2 e^{3}}
	+
	\int_{c_2e^3}^{(c_2e^3)\vee (1/t)} \d x\, \frac{2\,(x-c_2e^{3})}{x(\log x/c_2)^3}\ .
\end{align}
On the right-hand side, bound the terms by $t^{-1}(\log((e^3)\vee (1/c_2t)))^{-2}$ and $\frac{ 2 }{\log e^{3}} I'_2$, respectively.
Doing so gives $ I'_2\leq c\, e^{t}\, t^{-1} |\log (t\wedge\frac{1}{2})|^{-2}$.
As for $I_3$, since $x \leq 1/\e^{2}c_0$, the denominator of the integrand is $|\log(x/c_2)|^2$.
This is $\geq |\log (t\wedge\frac{1}{2})|^2/c$ because $x \geq (c_2e^{3})\vee (1/t)$.
For $I_4$, the denominator is $(\log \e^2-c)^2$, which is $\geq|\log (t\wedge\frac{1}{2})|^2/c$ because $t\geq \e^{2}c_0$.
Bounding the denominators this way and releasing the range of integration to $x\in[0,\infty)$ give $ I_3+ I_4 \leq c\, t^{-1} |\log (t\wedge\frac{1}{2})|^{-2}$.

\providecommand{\bysame}{\leavevmode\hbox to3em{\hrulefill}\thinspace}
\providecommand{\MR}{\relax\ifhmode\unskip\space\fi MR }
\providecommand{\MRhref}[2]{%
	\href{http://www.ams.org/mathscinet-getitem?mr=#1}{#2}
}
\providecommand{\href}[2]{#2}

\begin{acks}
 The authors would like to thank two anonymous referees for their careful reading and comments which improved the paper.
\end{acks}

\end{document}